\documentclass[reqno]{amsproc}
\usepackage[utf8]{inputenc}
\usepackage{cite}
\usepackage{xcolor}
\usepackage{color}
\usepackage{calligra}
\usepackage[T1]{fontenc}
\usepackage{latexsym}
\usepackage[all]{xy}
\usepackage[normalem]{ulem}
\usepackage{hyperref}
\hypersetup{
    colorlinks=true,		% false: boxed links; true: colored links
  linkcolor=blue,		% color of internal links
  citecolor=blue,		% color of links to bibliography
  filecolor=magenta,		% color of file links
  urlcolor=blue,		% color of internet links
  bookmarksdepth=4
}
\usepackage{mathrsfs}
\usepackage{amsmath,amsfonts,amssymb,amsxtra}
\numberwithin{equation}{section}
\newtheorem{theorem}{\bf Theorem}[section]

\newtheorem{lem}{\bf Lemma}[section]
\newtheorem{cor}{\bf Corollary}[section]
\newtheorem{remark}{\bf Remark}[section]

\newtheorem{example}{\bf Example}[section]

\newcommand\dv{\mathrm{div}}
\newcommand\tr{\mathrm{tr}}
\newcommand\rc{\mathrm{Ric}}
\usepackage{algorithmic}
\begin{document}

\title[Inequalities for eigenvalues of operators in divergence form]{Inequalities for eigenvalues of operators in divergence form on Riemannian manifolds isometrically immersed in Euclidean space}

\author{Cristiano S. Silva$^1$}
\author{Juliana F.R. Miranda$^2$}
\author{Marcio C. Ara\'ujo Filho$^3$}
\address{$^{1,2}$Departamento de Matem\'atica, Universidade Federal do Amazonas, Av. General Rodrigo Oct\'avio, 6200, 69080-900 Manaus, Amazonas, Brazil.}
\address{$^3$Departamento de Matemática, Universidade Federal de Rondônia, Campus Ji-Paraná, R. Rio Amazonas, 351, Jardim dos Migrantes, 76900-726 Ji-Paraná, Rondônia, Brazil}

\email{$^1$cristianosilva@ufam.edu.br}
\email{$^2$jfrmiranda@ufam.edu.br}
\email{$^3$marcio.araujo@unir.br}

\urladdr{$^{1,2}$http://dmice.ufam.edu.br}
\urladdr{$^3$https://dmejp.unir.br}


\begin{abstract}
In this paper, we compute universal inequalities of eigenvalues of a large class of second-order elliptic differential operators in divergence form,  that includes, e.g., the Laplace and Cheng-Yau operators, on a bounded domain in a complete Riemannian manifolds  isometrically immersed in Euclidean space. A key step in order to obtain the sequence of our estimates is to get the right Yang-type first inequality. We also prove some inequalities for manifolds supporting some special functions and tensors.
\end{abstract}
\maketitle

\section{Introduction}

Let $(M^n,\langle,\rangle)$ be an $n$-dimensional complete Riemannian manifold, and $\Omega~\subset~M^n$ be a bounded domain with smooth boundary $\partial\Omega$. Let us consider a function $\eta\in C^2(M)$ and a symmetric positive definite $(1,1)$-tensor $T$ on $M^n$. Since $\Omega$ is a bounded domain, there exist two positive real constants $\varepsilon$ and $\delta$, such that $\varepsilon \leq\langle T(X), X \rangle\leq\delta$, for any unit vector field $X$ on $\Omega$.

In this paper, we obtain some universal inequalities for eigenvalues of the following eigenvalue problem with the Dirichlet boundary condition:
\begin{equation}\label{problem1}
    \left\{\begin{array}{ccccc}
    \mathcal{L}  u  &=& - \lambda u & \mbox{in } & \Omega,\\
     u&=&0 & \mbox{on} & \partial\Omega,
    \end{array}
    \right.
\end{equation}
where $\mathcal{L}$ is defined as the second-order elliptic differential operator in the $(\eta,T)$-divergence form
\begin{equation}\label{eq1.1}
    \mathcal{L}u =\dv_\eta (T(\nabla u)) = \dv(T(\nabla u)) - \langle \nabla \eta, T(\nabla u) \rangle.
\end{equation}
Here $\dv$ stands for the divergence of smooth vector fields and $\nabla$ for the gradient of smooth functions. We highlight that the operator $\mathcal{L}$  is a class of second-order elliptic differential operators in divergence form that includes, e.g., the Laplace-Beltrami and Cheng-Yau operators with an approach using tensors.

We observe that $\mathcal{L}$ is a formally self-adjoint operator in the Hilbert space $\mathcal{H}_0^1(\Omega,e^{-\eta}d\Omega)$ of all functions in $L^2(\Omega,e^{-\eta}d\Omega)$ that vanish on $\partial \Omega$ in the sense of the trace, see Section~\ref{preliminaries}. Thus, Problem~\eqref{problem1} has a real and discrete spectrum
\begin{equation}\label{spectrum}
    0 < \lambda_{1} \leq \lambda_{2} \leq \cdots \leq \lambda_{k} \leq \cdots\to\infty,
\end{equation}
where each $\lambda_{i}$ is repeated according to its multiplicity. Eigenspaces belonging to distinct eigenvalues are orthogonal in $L^2(\Omega,e^{-\eta}d\Omega)$, which is the direct sum of all the eigenspaces. We refer to the dimension of each eigenspace as the multiplicity of the eigenvalue, for more details see Chavel~\cite{chavel}. 

Problem~\eqref{problem1} is an extension of the so-called fixed membrane problem, which in the case of two-dimensional euclidean space its eigenvalues are related to the characteristic vibrational frequencies of a uniformly stretched homogeneous membrane in the shape of  $\Omega$ with a fixed frontier, see Problem~\eqref{problem-3.4} for $\Omega\subset \mathbb{R}^2$. Problems relating the shape of a domain to the spectrum of an elliptic differential operator in divergence form are among the most fascinating of mathematical analysis. In this sense, one of the reasons which make them so interesting is that they related diﬀerent ﬁelds of mathematics such as Spectral Theory, Riemannian Geometry, and Partial Diﬀerential Equations. Furthermore, it is important to resalt that operators in divergence forms may play a keystone role in the understanding of countless physical facts.

It is natural and important to obtain universal inequalities for eigenvalues of the problems which involve elliptic differential operators in divergence form with  Dirichlet boundary condition on a bounded domain in a complete Riemannian manifold, at first, these inequalities are called universal because they do not involve domain dependencies. For these problems many mathematicians have been studying universal inequalities for eigenvalues in Euclidean and non-Euclidean environments, e.g.,  \cite{ChenCheng08}, \cite{ChengYang05,ChengYang06, ChengYang02, ChengYang09, ChengYang}, \cite{doCarmoWangXia}, \cite{HarrellStubbe}, \cite{LevitinParnovski}, \cite{XiaXu}, \cite{Yang} and their references. However, in the case of a complete Riemannian manifold is very difficult to find an appropriate trial function with ``nice'' properties such that one can infer universal inequalities for eigenvalues for the referred problem. Then we shall make use of the well-known Nash's theorem to construct appropriate trial functions with convenient properties and to obtain our results. In this regard, our results depend on a certain type of the mean curvature of the immersion, see Theorem~\ref{theorem_1.1}.

For the case of $T$ be the identity operator, and $\eta$ be a constant function, Problem~\eqref{problem1} becomes
\begin{equation}\label{problem-3.4}
    \left\{\begin{array}{ccccc}
    \Delta  u  &=& - \lambda u & \mbox{in } & \Omega,\\
     u&=&0 & \mbox{on} & \partial\Omega,
    \end{array}
    \right.
\end{equation}
 where $\Delta$ is the Laplace-Beltrami operator. Some interesting eigenvalue estimates for Problem~\eqref{problem-3.4} had been obtained, for instance, in Chen and Cheng~\cite{ChenCheng08}, Cheng and Yang~\cite{ChengYang02,ChengYang05,ChengYang06}, and Yang~\cite{Yang}, some of them will be referred in due course. We also emphasize here that there exists a beautiful and  abstract method used, for example, by Harrell II and Stubbe~\cite{HarrellStubbe}, Levitin and Parnovski~\cite{LevitinParnovski} and others to obtain eigenvalue estimates for Problem~\eqref{problem-3.4}, but we are not yet in a favorable condition to apply it in this case.

For the case of $T$ be the identity operator and $\eta$ be not necessarily a constant function, Problem~\eqref{problem1} becomes
\begin{equation}\label{problem-3.5}
    \left\{\begin{array}{ccccc}
    \Delta_\eta  u  &=& - \lambda u & \mbox{in } & \Omega,\\
     u&=&0 & \mbox{on} & \partial\Omega,
    \end{array}
    \right.
\end{equation}
where $\Delta_\eta$ is the drifted Laplace-Beltrami operator. In this case, it is worth mentioning the work by the second author and Gomes~\cite{GomesMiranda}, and the paper by Xia and Xu~\cite{XiaXu} in which we can find some estimates of the eigenvalues of Problem~\eqref{problem-3.5}.
Later, we can see that our result generalized some eigenvalues estimates of Problem~\eqref{problem-3.5} proved in these papers, see Remark~\ref{remark_1}.

We will call $T$ a {\it divergence-free} tensor when $\dv{T} = 0$. Divergence-free tensors often arise from physical facts, and the interested reader can find some of them on Serre~\cite{Serre}. For this case, the operator $\mathcal{L}$ becomes $\mathcal{L}f = \square f - \langle\nabla \eta, T(\nabla f)\rangle$, where $\square$ is the Cheng-Yau operator~\cite{ChengYau}, this shows that $\mathcal{L}$ is a first-order perturbation of the  Cheng-Yau operator $\square$, i.e., $\mathcal{L}$ is a  \emph{drifted Cheng-Yau operator} with a \emph{drifting function} $\eta$, which we will denote by $\square_\eta$, for more details see Araújo Filho and Gomes~\cite{AF-G} and \cite{AG}. For the case of $T$ be the divergence-free tensor and $\eta$ be not necessarily a constant function, Problem~\eqref{problem1} becomes
\begin{equation}\label{problem_CY}
    \left\{\begin{array}{ccccc}
    \square_\eta  u  &=& - \lambda u & \mbox{in } & \Omega,\\
     u&=&0 & \mbox{on} & \partial\Omega,
    \end{array}
    \right.
\end{equation}
where $\square_\eta$ is the drifted Cheng-Yau operator. Furthermore, in~\cite{AG} the authors have studied universal inequalities for eigenvalues of Problem~\eqref{problem_CY} and proved
\begin{equation}\label{drifted_CY_Ineq}
    \sum_{i=1}^k(\lambda_{k+1}-\lambda_{i})^2 \leq \frac{4\delta}{n\varepsilon}\sum_{i=1}^k(\lambda_{k+1}-\lambda_{i})\left(\lambda_{i}+\frac{n^2H_0^2+4C_0}{4\delta}\right),
\end{equation}
where $C_0$ is given by \eqref{C_0} and $H_0=\sup_\Omega |{\bf H}_T|$. Here and in what follows, the smooth vector field ${\bf H}_T$ on $M^n$ stands for the generalized mean curvature vector  associated with the $(1,1)$-tensor $T$, as defined in \eqref{HT}. 

Our first result is a quadratic estimate for the eigenvalues of Problem~\eqref{problem1} and is the fundamental key to get some of our bound of eigenvalues.

\begin{theorem}\label{theorem_1.1}
Let $\Omega$ be a bounded domain in an $n$-dimensional complete Riemannian manifold $M^n$ isometrically immersed in $\mathbb{R}^m$, and $\lambda_{i}$ be the $i$-th eigenvalue of Problem~\eqref{problem1}. Then, we have
\begin{equation*}
    \sum_{i=1}^k(\lambda_{k+1}-\lambda_{i})^2 \leq \frac{4\delta}{n\varepsilon}\sum_{i=1}^k(\lambda_{k+1}-\lambda_{i})\left(\lambda_{i} + \frac{n^2H_0^2+4C_0+T_0^2}{4\delta} \right),
\end{equation*}
where
\begin{equation}\label{C_0}
    C_0=\sup_\Omega \left\{\frac{1}{2}\dv \Big( T\big(T(\nabla \eta)-\tr(\nabla T)\big)\Big) - \frac{1}{4}|T(\nabla \eta)|^2\right\},
\end{equation}
$T_0=\sup_{\Omega}|\tr(\nabla T)|$, $\eta_0=\sup_{\Omega}|\nabla \eta|$ and $H_0=\sup_\Omega |{\bf H}_T|$. 
\end{theorem}

The quadratic inequality in Theorem~\ref{theorem_1.1} was identified as the keystone estimate for the applications of our results. 

\begin{remark}\label{remark_1}
Theorem~\ref{theorem_1.1} is an generalization of Inequality~\eqref{drifted_CY_Ineq}, since when $T$ is a divergence-free tensor we have $T_0 = 0$. Moreover, 
 Theorem~\ref{theorem_1.1} is an extension of Theorem~4 in \cite{GomesMiranda} which, in turn, is a slight modification from Theorem~1.2 in \cite{XiaXu}. We can see this fact by taking $T$ to be the identity operator, then  our result becomes an improvement of Theorem~4 in \cite{GomesMiranda}, since $4C_0=\sup_\Omega\Big\{2\Delta \eta-|\nabla \eta|^2\Big\}\leq 2\Bar{\eta_0}+\eta_0^2$, where $\Bar{\eta_0}=\sup_\Omega |\Delta_\eta \eta|$ and $\eta_0 = \sup_\Omega |\nabla \eta|$.
\end{remark}

We also give a general inequality for the lower-order eigenvalues of the second-order elliptic differential operator $\mathcal{L}$ in the $(\eta,T)$-divergence form defined in \eqref{eq1.1} as follows.

\begin{theorem}\label{theorem_1.2}
Let $\Omega$ be a bounded domain in an $n$-dimensional complete Riemannian manifold $M^n$ isometrically immersed in $\mathbb{R}^m$ and denote by $\lambda_{i}$ the $i$-th eigenvalue of Problem~\eqref{problem1}, for $i=1,\ldots,n$. Then, we have
\begin{equation*}
    \sum_{i=1}^n(\lambda_{i+1} - \lambda_{1}) \leq \frac{4\delta}{\varepsilon}\left(\lambda_{1} + \frac{n^2H_0^2+4C_0+T_0^2}{4\delta} \right),
\end{equation*}
 where the constants  $T_0, H_0$ and $C_0$ are as in Theorem~\ref{theorem_1.1}.
\end{theorem}
\begin{remark}
Theorems~\ref{theorem_1.1} and \ref{theorem_1.2} (see Corollary~\ref{cor-inf}) generalize two results by Chen and Cheng \cite[Theorem~1.1]{ChenCheng08} and ~\cite[Theorem~1.2]{ChenCheng08}, respectively.
\end{remark}

For applications of ours theorems, we set
\begin{equation*}
    \upsilon_{i}=\lambda_{i}+\frac{n^2H_0^2+4C_0+T_0^2}{4\delta}.
\end{equation*}
With this new notation, inequality in Theorem~\ref{theorem_1.1} can be rewritten as
\begin{equation}\label{T-parallel}
    \sum_{i=1}^k(\upsilon_{k+1}-\upsilon_{i})^2 \leq \frac{4\delta}{n\varepsilon}\sum_{i=1}^k(\upsilon_{k+1}-\upsilon_{i})\upsilon_{i}.
\end{equation}
Moreover, from \eqref{spectrum} we can see that $0<\upsilon_1 \leq \upsilon_2 \leq \cdots \leq \upsilon_{k}\leq \cdots\to\infty$. So, by applying the recursion formula of Cheng and Yang~\cite{ChengYang02} (see also \cite{ChengYang09} for most general case), we immediately obtain: 
\begin{cor}\label{recurs-ineq}
Under the same setup as in Theorem~\ref{theorem_1.1}, we have
\begin{equation*}
    \upsilon_{k+1} \leq \left(1+\frac{4\delta}{n\varepsilon}\right)k^{\frac{2\delta}{n\varepsilon}}\upsilon_1.
\end{equation*}
\end{cor}
When the $\mathcal{L}$ operator becomes the Laplace-Beltrami operator, we have $\varepsilon=\delta=1$, then from the classical Weyl’s asymptotic formula for the eigenvalues~\cite{Weyl}, we observe that the previous estimate is better bound in the sense of the order of $k$. The interested reader also can consult Ashabaugh's discussion~\cite[p. 15]{Ashbaugh} about this fact.

From quadratic inequality \eqref{T-parallel} and following the steps of the proof of Theorem~4 in Gomes and Miranda~\cite{GomesMiranda},  we also get.
\begin{cor}\label{cor-yang-type-ineq}
Under the same setup as in Theorem~\ref{theorem_1.1}, we have
\begin{align*}
 \upsilon_{k+1}\leq\frac{1}{k}\left(1 + \frac{4\delta}{n\varepsilon}		\right)\sum_{i=1}^k\upsilon_{i}+\left[\Big( \frac{2\delta}{kn\varepsilon}\sum_{i=1}^k\upsilon_{i}\Big)^{2} - \frac{1}{k}\Big(1 + \frac{4\delta}{n\varepsilon} \Big)\sum_{j=1}^k\Big(\upsilon_j-\frac{1}{k}\sum_{i=1}^k\upsilon_{i}\Big)^{2} \ \right]^\frac{1}{2}
\end{align*}
and
\begin{equation*}
 \upsilon_{k+1}-\upsilon_{k}\leq2 \left[\left( \frac{2\delta}{kn\varepsilon}\sum_{i=1}^k\upsilon_{i}
\right)^{2} - \frac{1}{k}\left(1 + \frac{4\delta}{n\varepsilon} \right) \sum_{j=1}^k\left(\upsilon_j-\frac{1}{k}\sum_{i=1}^k\upsilon_{i}\right)^{2}\right]^\frac{1}{2}.
\end{equation*}
\end{cor}

Furthermore, with this notation, Theorem~\ref{theorem_1.2} can be rewritten as:
\begin{cor}\label{cor-inf}
Under the same setup as in Theorem~\ref{theorem_1.2}, we have
\begin{equation*}
    \frac{\upsilon_2+\upsilon_3+\cdots+\upsilon_{n+1}}{\upsilon_1}\leq n+\frac{4\delta}{\varepsilon}.
\end{equation*}
\end{cor}

\begin{remark}
Corollaries~\ref{recurs-ineq},~\ref{cor-yang-type-ineq} and ~\ref{cor-inf}  generalize \cite[Corollaries~1.2,~1.3 and 1.5]{AG}  by third author and Gomes, respectively. 
\end{remark}

Also, under the same setup as in Theorem~\ref{theorem_1.1}, since $\upsilon_{i}\leq \upsilon_{i+1}$ for all $i= 1, 2, \ldots$ we know that $1/k\Big(\sum_{i=1}^k\upsilon_i \Big)\leq \sum_{i=1}^k \upsilon_k^2 $ and from \eqref{T-parallel} we can get a generalization of the Yang-type second inequality, as follows:
\begin{align*}
   \upsilon_{k+1}\leq\frac{1}{k}\left(1 + \frac{4\delta}{n\varepsilon}
			\right)\sum_{i=1}^k\upsilon_i.
\end{align*} 
This inequality generalizes Chen and Cheng~\cite[Corollary~1.1]{ChenCheng08}.

To conclude this section, we will present the following result of eigenvalues estimates for Problem~\eqref{problem1} in Riemannian manifolds admitting special functions. Such types of results were based on the Comparison theorem and initially obtained by do Carmo et al.~\cite{doCarmoWangXia}.
\begin{theorem}\label{theorem3}
Let $M$ be an $n$-dimensional complete Riemannian manifold and let $\Omega$ be a bounded domain with smooth boundary in $M$.  We denote by $\lambda_{i}$ the $i$-th eigenvalue of Problem~\eqref{problem1}.
\begin{enumerate}
    \item[i)] If there exists a function $\theta: \Omega \to \mathbb{R}$ and a constant $A_0$ such that
\begin{equation}\label{specialfunction-1}
    |\nabla \theta| = 1 \quad \mbox{and} \quad |\Delta\theta| \leq A_0, \quad \mbox{on} \quad \Omega,
\end{equation}
with $\nabla \theta$ an eigenfunction of $T$, then
\begin{equation*}
    \sum_{i=1}^k(\lambda_{k+1}-\lambda_{i})^2 \leq \frac{8\delta}{\varepsilon}\sum_{i=1}^k(\lambda_{k+1}-\lambda_{i})\Big(\lambda_{i}+ \frac{\delta(A_0+\eta_0)^2}{4}\Big).
\end{equation*}
\item[ii)] If there exists a function $\psi: \Omega \to \mathbb{R}$ and a constant $B_0$ such that
\begin{equation}\label{specialfunction-3}
    |\nabla \psi| = 1 \quad \mbox{and} \quad \mathcal{L}\psi = B_0, \quad \mbox{on} \quad \Omega,
\end{equation}
then
\begin{equation*}
    \sum_{i=1}^k(\lambda_{k+1}-\lambda_{i})^2 \leq \frac{4\delta}{\varepsilon}\sum_{i=1}^k(\lambda_{k+1}-\lambda_{i})\Big(\lambda_{i}- \frac{B_0^2}{4}\Big).
\end{equation*}
In particular, the first eigenvalue satisfies $\lambda_{1}(\Omega ) \geq \frac{B_0^2}{4}$.
\item[iii)] If $M$ admits an map $f=(f_1, f_2, \ldots, f_{m+1}): \Omega \to \mathbb{S}^m$ such that
\begin{equation}\label{specialfunction-2}
    \sum_{\ell =1}^{m+1}f_{\ell}^2 = 1 \quad \mbox{and} \quad \mathcal{L} f_{\ell} = - \gamma f_{\ell}, \quad  \mbox{for} \quad \ell = 1, \ldots, m+1,
\end{equation}
for some constant $\gamma$, then
\begin{equation*}
    \sum_{i=1}^k(\lambda_{k+1}-\lambda_{i})^2 \leq \frac{4\delta}{\varepsilon}\sum_{i=1}^k(\lambda_{k+1}-\lambda_{i})\Big(\lambda_{i} + \frac{\gamma\varepsilon}{4\delta}\Big),
\end{equation*}
where $\mathbb{S}^m$ is the unit $m$-sphere.
\end{enumerate}
\end{theorem}

\begin{remark}
The function $\theta$ in item i) and $f_1, . . . , f_{m+1}$ in item iii) of Theorem~\ref{theorem3} are only required to be defined on $\Omega$ and the corresponding inequalities for the eigenvalues depend on $\Omega$. On the other hand, for those manifolds on which there exist globally defined functions of similar nature, the inequalities for eigenvalues of $\mathcal{L}$ are independent of the domains and, therefore, are universal.
\end{remark}

In what follows we give examples of special functions satisfying the conditions in Theorem~\ref{theorem3}. For deeper details, the interested reader can consult \cite{doCarmoWangXia}, \cite{WangXia}, \cite{XiaXu} and their references.
\begin{example}
Let $N^{n-1}$ be a complete Riemannian manifold with Ricci non-negative. The warped product $M^n=\mathbb{R}\times_{e^t}N^{n-1}$ with the warped product metric $ds_M^2 = dt^2 + e^{2t}ds_N^2$ is a complete Riemannian manifold and $\rc_M \geq -(n-1)$. Therefore, the function $\theta: M \to \mathbb{R}$ given by $\theta(p, t) = t$ satisfies $|\nabla \theta|=1$ and $|\Delta \theta| \leq n-1$, that is, satisfy the conditions in item i) of Theorem~\ref{theorem3}.
\end{example}
\begin{example}\label{example-2}
Let $\mathbb{H}^n(-1)$ the hyperbolic space with the upper half-plane model, that is, $\mathbb{H}^n(-1) = \{x=(x_1, \ldots, x_n)\in\mathbb{R}^n; x_n > 0\}$  with the standard metric $g_{ij}=x_n^{-2}\delta_{ij}$. The function $f(x)=\ln x_n$ satisfy
$|\nabla f|=1$ and $\Delta f = 1 -n$, that is, satisfy the conditions in item ii) when $T=I$. Moreover, in the most general case, if the drifting function $\eta$ is radially constant and $T( \nabla \ln x_n) = \psi  \nabla \ln x_n$  for some radially constant function $\psi \in C^\infty(\Omega)$ on a bounded domain $\Omega \in \mathbb{H}^n(-1)$, we can get $\mathcal{L} f = (1-n)\psi$ which satisfy the conditions in item ii) of Theorem~\ref{theorem3} when $\psi$ is a constant function.
\end{example}
\begin{example}
For any compact homogeneous Riemannian manifold we can give eigenmaps to some unit sphere for the first positive eigenvalue of the Laplacian and so satisfy the conditions in item iii) of Theorem~\ref{theorem3} when $T=I$.
\end{example}
Immediately from item ii) of Theorem~\ref{theorem3} and Example~\ref{example-2} we obtain. 
\begin{cor}
For a bounded domain $\Omega$ in $\mathbb{H}^n(-1)$, eigenvalues $\lambda_{i}$'s of Laplacian with Dirichlet boundary condition satisfy
\begin{equation*}
    \sum_{i=1}^k(\lambda_{k+1}-\lambda_{i})^2 \leq 4\sum_{i=1}^k(\lambda_{k+1}-\lambda_{i})\Big(\lambda_{i}- \frac{(n-1)^2}{4}\Big).
\end{equation*}
In particular, $\lambda_{1} \geq \frac{(n-1)^2}{4}$.
\end{cor}
The universal inequality above was obtained by Cheng and Yang in \cite[Theorem~1.2]{ChengYang09}. 

\section{Preliminaries}\label{preliminaries}

In this section, we establish some basic notations and describe properties of a $(1,1)$-tensor in a bounded domain $\Omega\subset M^n$ with smooth boundary $\partial\Omega$.

Throughout this work, we are frequently using the identification of a $(0,2)$-tensor $T:\mathfrak{X}(M)\times\mathfrak{X}(M)\to~C^{\infty}(M)$ with its associated $(1,1)$-tensor $T:\mathfrak{X}(M)\to\mathfrak{X}(M)$ by the equation
\begin{equation*}
    \langle T(X), Y \rangle = T(X, Y).
\end{equation*}
Moreover, the tensor $\langle , \rangle$ will be identified with the identity $I$ in $\mathfrak{X}(M)$. Since, 
$\varepsilon \leq\langle T(X), X \rangle\leq\delta$, for any unit vector field $X$ on $\Omega$, we have
\begin{align}\label{T-property}
\varepsilon  \langle T(Y), Y\rangle \leq |T(Y)|^2 \leq \delta  \langle T(Y), Y\rangle \quad \mbox{for all} \quad Y \in \mathfrak{X}(M).
\end{align}
Consequently, 
\begin{align*}
    \varepsilon^2|\nabla \eta|^2 \leq |T(\nabla\eta)|^2 \leq \delta^2|\nabla\eta|^2.
\end{align*}

For an $n$-dimensional complete Riemannian manifold $(M^n, \langle , \rangle)$  isometrically immersed in $\mathbb{R}^m$ we denote by $\alpha$ its second fundamental form and by ${\bf H}=\frac{1}{n}\tr (\alpha)$ its mean curvature vector. Let $\{e_1, e_2, \ldots, e_n\}$ be an orthonormal basis of $T_pM$, for a symmetric $(1, 1)$-tensor $T$, we define the generalized mean curvature vector at $p\in M$ as 
\begin{equation}\label{HT}
    {\bf H}_T=\frac{1}{n}\sum_{i,j=1}^nT(e_{i}, e_{j})\alpha(e_{i}, e_{j})=\frac{1}{n}\sum_{i=1}^n\alpha(T(e_{i}), e_{i})=\frac{1}{n}\tr{(\alpha\circ T)}.
\end{equation}
 We shall use the notations
\begin{align*}
    \tr(\nabla T):=\sum_{i=1}^n(\nabla T)(e_{i}, e_{i}), \quad |T|=\Big(\sum_{i=1}^n|T(e_{i})|^2\Big)^{\frac{1}{2}}
\end{align*}
 and $\|\cdot\|_{L^2}$ for the canonical norm of a function in $L^2(\Omega,e^{-\eta}d\Omega)$. 

It is important to observe that, from definition of $\eta$-divergence of $X$ (see Eq.~\eqref{eq1.1}), we have
\begin{equation*}
    \dv_\eta(fX)=f\dv_\eta X + \langle \nabla f, X \rangle,
\end{equation*}
and then
\begin{align}\label{property1}
    \mathcal{L}(f\ell) = f\mathcal{L}\ell + 2T(\nabla f, \nabla \ell) + \ell\mathcal{L}f
\end{align}
for all $f,\ell \in C^\infty(M)$. Moreover, the divergence theorem is valid as follows:
\begin{equation*}
   \int_\Omega\dv_\eta X\mathrm{dm} = \int_{\partial\Omega} \langle X, \nu\rangle d\mu,
\end{equation*}
where $\mathrm{dm} = e^{-\eta}d\Omega$ is the weighted volume form on $\Omega$ and $d\mu = e^{-\eta}d\partial\Omega$ is the weighted area form on $\partial\Omega$ induced by the 
outward pointing unit normal vector field $\nu$ along $\partial \Omega$. In particular, by taking $X=T(\nabla f)$, we get 
\begin{equation*}
   \int_\Omega\mathcal{L}{f}\mathrm{dm} = \int_{\partial\Omega}T(\nabla f, \nu) d\mu
\end{equation*}
and the integration by parts formula:
\begin{equation}\label{parts}
     \int_{\Omega}\ell\mathcal{L}{f}\mathrm{dm} =-\int_\Omega T(\nabla\ell, \nabla f)\mathrm{dm} + \int_{\partial\Omega}\ell T(\nabla f, \nu) d\mu
\end{equation}
for all $\ell, f \in C^\infty(M)$. Therefore, from \eqref{parts} we can see that $\mathcal{L}$ is a formally self-adjoint operator in the Hilbert space $\mathcal{H}_0^1(\Omega,\mathrm{dm})$.

For the eigenfunction $u_{i}$ corresponding to the eigenvalue $\lambda_{i}$, from Problem~\eqref{problem1} and integration by parts formula~\eqref{parts}, we have
\begin{equation*}
    \lambda_{i}\int_\Omega u_{i}^{2}\mathrm{dm} = -\int_\Omega u_{i}\mathcal{L}u_{i}\mathrm{dm} = \int_\Omega T(\nabla u_{i}, \nabla u_{i})\mathrm{dm}.
\end{equation*}

Finally, we should mention a significant work for us with respect to Problem~\eqref{problem1} by the second author and Gomes~\cite[Section~2]{GomesMiranda} where is possible to find some geometric motivations to work with the operator $\mathcal{L}$ in the $(\eta,T)$-divergence form  and a Bochner-type formula for it on Riemannian manifolds.

\section{Proof of Theorems~\ref{theorem_1.1}, \ref{theorem_1.2} and \ref{theorem3}}

In order to prove our results, we will need two keystone technical lemmas. The first, was proved by Gomes and Miranda~\cite[Proposition~1]{GomesMiranda}.
\begin{lem}[Gomes and Miranda~\cite{GomesMiranda}, p. 7-9]\label{lemma1}
Let $\Omega$ be a bounded domain in an $n$-dimensional complete Riemannian manifold $M$. Let $\lambda_{i}$ be the i-th eigenvalue of Problem~\eqref{problem1} and $u_{i}$ its correspoding normalized real-valued eigenfunction. Then, for any $f\in C^2(\Omega)\cap C^1(\partial \Omega)$, we obtain
\begin{align*}
    \sum_{i=1}^k(\lambda_{k+1}-\lambda_{i})^2 \int_{\Omega}T(\nabla f, \nabla f)u_{i}^{2}\mathrm{dm} \leq 4\sum_{i=1}^k(\lambda_{k+1}-\lambda_{i})\|T(\nabla f, \nabla u_{i}) + \frac{1}{2}u_{i}\mathcal{L}f \|^2_{L^2}.
\end{align*}
\end{lem}

The second lemma was proved by Araújo Filho and Gomes~\cite{AG}  using Lemma~\ref{lemma1} and following the steps of the proof of Proposition~2 in Gomes and Miranda~\cite{GomesMiranda} with appropriate adaptations.

\begin{lem}[Araújo Filho and Gomes~\cite{AG}, p. 8-9]\label{lemma2}
Let $\Omega$ be a bounded domain in an $n$-dimensional complete Riemannian manifold $M$ isometrically immersed in $\mathbb{R}^m$, $\lambda_{i}$ be the $i$-th eigenvalue of Problem~\eqref{problem1} and $u_{i}$ its corresponding $L^2(\Omega,\mathrm{dm})$-normalized real-valued eigenfunction. Then is valid
\begin{align*}
   &\sum_{i=1}^k (\lambda_{k+1}-\lambda_{i})^2 \int_\Omega \tr(T)u_{i}^{2}\mathrm{dm}\\
      &\leq 4\sum_{i=1}^k(\lambda_{k+1}-\lambda_{i})\Bigg\{\|T (\nabla u_{i})\|^2_{L^2} +\frac{n^2}{4}\int_{\Omega}u_{i}^{2}|{\bf H}_T|^2\mathrm{dm}+ \int_\Omega u_{i} T(\tr(\nabla T), \nabla u_{i})\mathrm{dm}\\
   &+\int_\Omega\!\!u_{i}^{2}\Big(\frac{1}{2}\dv (T^2(\nabla \eta))\!-\!\frac{1}{4}|T(\nabla \eta)|^2\Big)\mathrm{dm}+\frac{1}{4}\!\int_\Omega u_{i}^{2}\langle\tr(\nabla T), \tr(\nabla T) - 2T(\nabla \eta)\rangle\mathrm{dm} \Bigg\}.
\end{align*}
\end{lem}
Now, we are in a position to give the proof of the theorems of this paper. 

\subsection{Proof of Theorem~\ref{theorem_1.1}}
\begin{proof}
	The proof is a consequence of Lemma~\ref{lemma2}. We initially observe that, since there are positive real numbers $\varepsilon$ and $\delta$ such that $\varepsilon|X|^{2}\leq \langle T(X), X \rangle \leq\delta|X|^{2} $, then for any vector field $ X\in \mathfrak{X}(M) $, we have $ n\varepsilon \leq \mathrm{tr}(T) $ and consequently
	\begin{align}\label{3,1}
		n\varepsilon = n\varepsilon\int_{\Omega}u_{i}^{2}\mathrm{dm} = \int_{\Omega}n\varepsilon u_{i}^{2}\mathrm{dm} \leq \int_{\Omega}\mathrm{tr}(T)u_{i}^{2}\mathrm{dm}.
	\end{align}
	Since we set $ \displaystyle H_{0}  = \sup_{\overline{\Omega}}|{\bf{H}}_{T}| $ and using \eqref{T-property}, we obtain
	\begin{align}
		\nonumber
		\int_{\Omega}|T(\nabla u_{i})|^{2} \mathrm{dm} + \frac{n^{2}}{4}\int_{\Omega}u_{i}^{2}|{\bf{H}}_{T}|^{2}\mathrm{dm} &\leq \int_{\Omega}\delta\left\langle T(\nabla u_{i}), \nabla u_{i} \right\rangle \mathrm{dm} + \frac{n^{2}}{4}\int_{\Omega}u_{i}^{2}H_{0}^{2}\mathrm{dm}\\
		\nonumber
		&\leq \delta\!\int_{\Omega}\!\left\langle T(\nabla u_{i}), \nabla u_{i} \right\rangle  \mathrm{dm} + \frac{n^{2}H_{0}^{2}}{4}\int_{\Omega}u_{i}^{2}\mathrm{dm}  \\
		&\leq \delta\lambda_{i} + \frac{n^{2}H_{0}^{2}}{4}.
	\end{align}
	Let us consider $ \displaystyle T_{0} = \sup_{\overline{\Omega}}|\mathrm{tr}(\nabla T)| $, to obtain 
	\begin{align}
		\nonumber
		\int_{\Omega}u_{i}&\left\langle \mathrm{tr}(\nabla T),T(\nabla u_{i})\right\rangle \mathrm{dm} + \frac{1}{4}\int_{\Omega}u_{i}^{2}\left\langle \mathrm{tr}(\nabla T), \mathrm{tr}(\nabla T) - 2T(\nabla \eta)\right\rangle\mathrm{dm} \\
		\nonumber
		=& \frac{1}{2}\int_{\Omega}\left\langle T(\mathrm{tr}(\nabla T)),\nabla (u_{i}^{2})\right\rangle \mathrm{dm} + \frac{1}{4}\int_{\Omega}u_{i}^{2}\left\langle \mathrm{tr}(\nabla T), \mathrm{tr}(\nabla T) - 2T(\nabla \eta)\right\rangle\mathrm{dm} \\
		\nonumber
		=& -\frac{1}{2}\int_{\Omega}u_{i}^{2} div_{\eta}(T(\mathrm{tr}(\nabla T))) \mathrm{dm} + \frac{1}{4}\int_{\Omega}u_{i}^{2} |\mathrm{tr}(\nabla T)|^{2} \mathrm{dm} \\ 
        \nonumber
        &-\frac{1}{2}\int_{\Omega}u_{i}^{2}\left\langle \mathrm{tr}(\nabla T), T(\nabla \eta)\right\rangle \mathrm{dm} \\
		\nonumber
		=& - \frac{1}{2}\int_{\Omega}u_{i}^{2} div(T(\mathrm{tr}(\nabla T))) \mathrm{dm} + \frac{1}{2}\int_{\Omega}u_{i}^{2} \left\langle \mathrm{tr}(\nabla T), T(\nabla \eta) \right\rangle \mathrm{dm} \\
		\nonumber
		&\quad +\frac{1}{4}\int_{\Omega}u_{i}^{2} |\mathrm{tr}(\nabla T)|^{2} \mathrm{dm} - \frac{1}{2}\int_{\Omega}u_{i}^{2}\left\langle \mathrm{tr}(\nabla T), T(\nabla \eta)\right\rangle \mathrm{dm} \\
		\leq& - \frac{1}{2}\int_{\Omega}u_{i}^{2} div(T(\mathrm{tr}(\nabla T))) \mathrm{dm} + \!\frac{T_{0}^{2}}{4}.
	\end{align}
	Moreover, notice that 
	\begin{align}\label{3,4}
		\nonumber
		&\int_{\Omega}u_{i}^{2} \left(\frac{1}{2}div(T^{2}(\nabla\eta)) - \frac{1}{4}|T(\nabla\eta)|^{2}\right) \mathrm{dm} - \frac{1}{2}\int_{\Omega} u_{i}^{2} div(T(\mathrm{tr}(\nabla T))) \mathrm{dm} \\
		=& \int_{\Omega}u_{i}^{2} \left(\frac{1}{2}div
		\bigg(T\Big(T(\nabla\eta) - \mathrm{tr}(\nabla T)\Big)\bigg) - \frac{1}{4}|T(\nabla\eta)|^{2}\right) \mathrm{dm} \leq C_{0}.
	\end{align}
	where $ \displaystyle C_{0} = \sup_{\overline{\Omega}}\left\{\frac{1}{2}div
	\bigg(T\Big(T(\nabla\eta) - \mathrm{tr}(\nabla T)\Big)\bigg) - \frac{1}{4}|T(\nabla\eta)|^{2} \right\} $.
    Therefore, substituting \eqref{3,1}-\eqref{3,4} into Lemma~\ref{lemma2} we obtain
	\begin{align*}
		n\varepsilon \sum_{i=1}^k(\lambda_{k+1}-\lambda_{i})^{2} \leq 4\sum_{i=1}^k(\lambda_{k+1}-\lambda_{i}) \Bigg(\delta\lambda_{i} + \frac{n^{2}H_{0}^{2}}{4} + C_{0} + \frac{T_{0}^{2}}{4}  \Bigg),
	\end{align*}
	that is,
	\begin{align*}
		\sum_{i=1}^k(\lambda_{k+1}-\lambda_{i})^{2} \leq \frac{4\delta}{n\varepsilon }\sum_{i=1}^k(\lambda_{k+1}-\lambda_{i}) \Bigg(\lambda_{i} + \frac{n^{2}H_{0}^{2} + 4C_{0} + T_{0}^{2}}{4\delta} \Bigg),
	\end{align*}
 which completes the proof of Theorem~\ref{theorem_1.1}.
\end{proof}

\subsection{Proof of Theorem~\ref{theorem_1.2}}
\begin{proof}
The proof follows the steps of the proof of \cite[Theorem~1.2]{AG}. Let $x=\left(x_{1}, \ldots, x_{m}\right)$ be the position vector of the immersion of $M$ in $\mathbb{R}^{m}$. Let us consider the matrix $D=\left(d_{i j}\right)_{m \times m}$ where
\begin{align*}
d_{ij}:=\int_{\Omega}x_{i}u_{1}u_{j+1} \mathrm{dm}.
\end{align*}
From the orthogonalization of Gram and Schmidt, we know that there exists an upper triangle matrix $R=\left(r_{i j}\right)_{m \times m}$ and an orthogonal matrix $S=\left(s_{i j}\right)_{m \times m}$ such that $R=S D$, namely
\begin{align*}
r_{ij}=\sum_{k=1}^{m}s_{ik}d_{kj}=\sum_{k=1}^{m}s_{ik}\int_{\Omega}x_{k}u_{1}u_{j+1}\mathrm{dm}=\int_{\Omega}\left(\sum_{k=1}^{m}s_{ik}x_{k}\right)u_{1}u_{j+1}\mathrm{dm}=0,
\end{align*}
for $1 \leq j<i \leq m$. By setting $y_{i}=\sum_{k=1}^{m} s_{i k} x_{k}$, we have
\begin{align*}
\int_{\Omega}y_{i}u_{1} u_{j+1}\mathrm{dm}=0\quad\text{ for }\quad 1\leq j<i\leq m.
\end{align*}
Let us denote $a_{i}=\int_{\Omega}y_{i}\left|u_{1}\right|^{2}\mathrm{dm}$ and consider the real-valued functions $w_{i}$ given by
\begin{align*}
w_{i}=\left(y_{i}-a_{i}\right)u_{1}
\end{align*}
so that
\begin{align*}
\left.w_{i}\right|_{\partial\Omega}=0\quad\text{ and }\quad\int_{\Omega}w_{i}u_{j+1}\mathrm{dm}=0,\quad\text{ for any }\quad j=1,\ldots, i-1.
\end{align*}
Then, from Rayleigh-Ritz inequality, we have for $1\leq i\leq m$
\begin{align*}
\lambda_{i+1}\left\|w_{i}\right\|_{L^{2}}^{2}\leq-\int_{\Omega}w_{i}\mathcal{L} w_{i}\mathrm{dm}.
\end{align*}
and from expressions (4.9)-(4.13) in \cite{AG}, we have
\begin{align}\label{3,5}
\sum_{i=1}^{m}(\lambda_{i+1}-\lambda_{1})\int_{\Omega}\left|u_{1}\right|^{2}T(\nabla y_{i}, \nabla y_{i})\mathrm{dm}\leq4\sum_{i=1}^{m}\left\|\frac{1}{2}u_{1}\mathcal{L} y_{i}+T(\nabla y_{i},\nabla u_{1})\right\|_{L^{2}}^{2}.
\end{align}
Hence, using the definition of $y_{i}$ and the fact that $S$ is an orthogonal matrix, we obtain, similarly to Eq. (3.17)-(3.24) in \cite{GomesMiranda}, that
\begin{align}\label{eqS1}
\left.
\begin{array}{c}
\displaystyle
\sum_{\ell=1}^{m}\left|\nabla y_{\ell}\right|^{2}=n, \\
\\
\displaystyle
\sum_{\ell=1}^{m}\left(\mathcal{L}y_{\ell}\right)^{2}=n^{2}\left|\mathbf{H}_{T}\right|^{2}+|\operatorname{tr}(\nabla T)-T(\nabla \eta)|^{2},\\
\\
\displaystyle
\sum_{\ell=1}^{m}\left|T(\nabla y_{\ell}, \nabla u_{1})\right|^{2}=\sum_{\ell=1}^{m}\left|T\left(e_{\ell}, \nabla u_{1}\right)\right|^{2}=\left|T(\nabla u_{1})\right|^{2},\\
\\
\displaystyle
\sum_{\ell=1}^{m}\mathcal{L}y_{\ell}\left\langle T(\nabla y_{\ell}), \nabla u_{1}\right\rangle
=\left\langle T(\mathrm{tr}(\nabla T), \nabla u_{1}\right\rangle-\left\langle T^{2}(\nabla\eta), \nabla u_{1}\right\rangle.
\end{array}
\right\}
\end{align}
Since there exist positive real numbers $\varepsilon$ and $\delta$ such that $\varepsilon\langle X, X\rangle \leq\langle T(X), X\rangle \leq$ $\delta\langle X, X\rangle$, for any vector field $X$ on $\Omega$, from \cite[Inequality~(4.18)]{AG} we have
\begin{align}\label{3,7}
\varepsilon \sum_{i=1}^{n}\left(\lambda_{i+1}-\lambda_{1}\right)\leq\sum_{i=1}^{m}\left(\lambda_{i+1}-\lambda_{1}\right) T(\nabla y_{i}, \nabla y_{i}).
\end{align}
Now, for the right side in \eqref{3,5}, using the equalities in \eqref{eqS1}, we get
\begin{align}\label{3,8}
\nonumber
4\sum_{i=1}^{m}&\left\|\frac{1}{2}u_{1}\mathcal{L} y_{i}+T(\nabla y_{i},\nabla u_{1})\right\|_{L^{2}}^{2}=4\sum_{i=1}^{m}\int_{\Omega}\left(\frac{1}{2}u_{1}\mathcal{L}y_{i}+T(\nabla y_{i},\nabla u_{1})\right)^{2}\mathrm{dm}\\
\nonumber
=&\ 4\sum_{i=1}^{m}\int_{\Omega}\left(\frac{1}{4}u_{1}^{2}\left(\mathcal{L}y_{i}\right)^{2}+\left|T(\nabla y_{i}, \nabla u_{1})\right|^{2}+u_{1}\mathcal{L}y_{i}\left\langle T(\nabla y_{i}),\nabla u_{1})\right\rangle\right)\mathrm{dm}\\
\nonumber
=&\ \int_{\Omega}u_{1}^{2}\left(n^{2}\left|\mathbf{H}_{T}\right|^{2}+|\operatorname{tr}(\nabla T)-T(\nabla \eta)|^{2}\right)\mathrm{dm}+4\int_{\Omega}\left|T(\nabla u_{1})\right|^{2}\mathrm{dm} \\
&+4\int_{\Omega}u_{1}\Big(\left\langle T(\mathrm{tr}(\nabla T), \nabla u_{1}\right\rangle-\left\langle T^{2}(\nabla\eta), \nabla u_{1}\right\rangle\Big)\mathrm{dm}.
\end{align}
Since $u_{1}|_{\partial\Omega}=0$ using the divergence theorem, we have
\begin{align}\label{3,9}
\nonumber
4\int_{\Omega}u_{1}&\Big(\left\langle T(\mathrm{tr}(\nabla T), \nabla u_{1}\right\rangle-\left\langle T^{2}(\nabla\eta), \nabla u_{1}\right\rangle\Big)\mathrm{dm}\\
\nonumber
=&2\int_{\Omega}\left\langle T(\mathrm{tr}(\nabla T), \nabla u_{1}^{2}\right\rangle\mathrm{dm}-2\int_{\Omega}\left\langle T^{2}(\nabla\eta), \nabla u_{1}^{2}\right\rangle\mathrm{dm}\\
\nonumber
=&-2\int_{\Omega}u_{1}^{2}\dv_{\eta}(T(\mathrm{tr}(\nabla T))\mathrm{dm}+2\int_{\Omega}u_{1}^{2}\dv_{\eta}(T^{2}(\nabla\eta))\mathrm{dm}\\
\nonumber
=&-2\int_{\Omega}u_{1}^{2}\Big(\dv(T(\mathrm{tr}(\nabla T))-\left\langle T(\mathrm{tr}(\nabla T),\nabla\eta\right\rangle\Big)\mathrm{dm}\\
&+2\int_{\Omega}u_{1}^{2}\Big(\dv(T^{2}(\nabla\eta))-\left|T(\nabla\eta)\right|^{2}\Big)\mathrm{dm}.
\end{align}
So, from \eqref{3,8} and \eqref{3,9}
\begin{align}\label{3,10}
4\sum_{i=1}^{m}&\left\|\frac{1}{2}u_{1}\mathcal{L} y_{i}+T(\nabla y_{i},\nabla u_{1})\right\|_{L^{2}}^{2}\nonumber\\
\nonumber
=&\ 4\Bigg\{\|T (\nabla u_{1})\|^2_{L^2} +\frac{n^2}{4}\int_{\Omega}u_{1}^{2}|{\bf H}_T|^2\mathrm{dm}+\frac{1}{4}\!\int_\Omega u_{1}^{2}\left|\tr(\nabla T)\right|^{2}\mathrm{dm}\\
\nonumber
&+\int_\Omega\!\!u_{1}^{2}\bigg[\frac{1}{2}\mathrm{div} \bigg(T\Big(T(\nabla\eta) - \mathrm{tr}(\nabla T)\Big)\bigg) - \frac{1}{4}|T(\nabla\eta)|^{2}\bigg]\mathrm{dm} \Bigg\}\\
\leq&\ 4\Bigg\{\delta\lambda_{1} +\frac{n^2H_{0}^{2}}{4}+\frac{T_{0}^{2}}{4}+C_{0}\Bigg\},
\end{align}
where $H_{0}$, $T_{0}$ e $C_{0}$ are as in Theorem~\ref{theorem_1.1}. Finally, from \eqref{3,7}, \eqref{3,10} and \eqref{3,5} we obtain
\begin{align*}
\sum_{i=1}^k(\lambda_{i+1}-\lambda_1)\leq \frac{4\delta}{\varepsilon}\Bigg(\lambda_{1} + \frac{n^{2}H_{0}^{2} + 4C_{0} + T_{0}^{2}}{4\delta}\Bigg).
\end{align*}
\end{proof}

\subsection{Proof of Theorem~\ref{theorem3}}
\begin{proof} i) Since $\nabla \theta$ is an eigenfunction of $T$, there are a constant $\gamma_0$ such that $T(\nabla \theta)=\gamma_0 \nabla \theta$ with $\gamma_0 \leq \delta$, then
\begin{equation*}
 \mathcal{L}\theta = \dv_\eta{(T(\nabla \theta))}=\gamma_0\dv_\eta{\nabla \theta}= \gamma_0\Delta_\eta \theta,   
\end{equation*}
and from \eqref{specialfunction-1} we have
\begin{equation*}
    |\mathcal{L}\theta| = |\gamma_0\Delta_\eta \theta|\leq \delta |\Delta_\eta \theta|=\delta|\Delta \theta - \langle \nabla \eta, \nabla \theta \rangle| \leq \delta (|\Delta \theta| +|\nabla \eta||\nabla \theta|) \leq \delta (A_0 + \eta_0).
\end{equation*}
Then, taking $f=\theta$ into Lemma~\ref{lemma1}, using the previous inequality and the Schwarz inequality, we obtain
\begin{align*}
    \varepsilon \sum_{i=1}^k&(\lambda_{k+1}-\lambda_{i})^2\\
    & \leq \sum_{i=1}^k(\lambda_{k+1}-\lambda_{i})^2 \int_\Omega T(\nabla \theta, \nabla \theta) u_{i}^{2}\mathrm{dm} \\
    &\leq \sum_{i=1}^k(\lambda_{k+1}-\lambda_{i}) \int_\Omega (u_{i} \mathcal{L}\theta + 2T(\nabla \theta, \nabla u_{i}))^2 \mathrm{dm} \\
    &\leq\sum_{i=1}^k(\lambda_{k+1}-\lambda_{i})\int_\Omega \Big(u_{i}^{2}(\mathcal{L}\theta)^2 + 4u_{i}\mathcal{L}\theta T(\nabla \theta, \nabla u_{i}) + 4(T(\nabla \theta, \nabla u_{i}))^2 \Big)\mathrm{dm}\\
    &\leq  \sum_{i=1}^k(\lambda_{k+1}-\lambda_{i}) \int_\Omega 2\Big(u_{i}^{2}(\mathcal{L}\theta)^2 + 4\langle \nabla \theta, T(\nabla u_{i})\rangle^2 \Big)\mathrm{dm}\\
   % &\leq 2 \sum_{i=1}^k(\lambda_{k+1}-\lambda_{i}) \int_\Omega \Big(u_{i}^{2}A_0^2 + 4\langle \nabla \theta, T(\nabla u_{i})\rangle^2 \Big)\mathrm{dm}\\
    &\leq 2\sum_{i=1}^k(\lambda_{k+1}-\lambda_{i}) \int_\Omega\Big(u_{i}^{2}\delta^2(A_0+\eta_0)^2 + 4|\nabla \theta|^2 |T(\nabla u_{i})|^2 \Big)\mathrm{dm}\\
     &= 8 \sum_{i=1}^k(\lambda_{k+1}-\lambda_{i})  \Big(\frac{\delta^2(A_0+\eta_0)^2}{4} + \|T(\nabla u_{i})\|_{L^2}^2 \Big)\\
      &\leq 8 \sum_{i=1}^k(\lambda_{k+1}-\lambda_{i})  \Big(\delta \lambda_{i} + \frac{\delta^2(A_0+\eta_0)^2}{4}\Big),
\end{align*}
which complete the proof of item i).\\
ii) Taking $f=\psi$  into Lemma~\ref{lemma1} we have
\begin{align*}
    \sum_{i=1}^k(\lambda_{k+1}-\lambda_{i})^2 \int_{\Omega}T(\nabla \psi, \nabla \psi)u_{i}^{2}\mathrm{dm} \leq 4\sum_{i=1}^k(\lambda_{k+1}-\lambda_{i})\|T(\nabla \psi, \nabla u_{i}) + \frac{1}{2}\mathcal{L}\psi u_{i}\|^2_{L^2}.
\end{align*}
From integration by parts, we get
\begin{align*}
     4\|T(\nabla \psi, \nabla u_{i}) + &\frac{1}{2}\mathcal{L}\psi u_{i}\|^2_{L^2}\\
     &= \int_\Omega (u_{i} \mathcal{L}\psi + 2T(\nabla \psi, \nabla u_{i}))^2 \mathrm{dm}\\
     &= \int_\Omega \Big(u_{i}^{2} (\mathcal{L}\psi)^2 +4u_{i}\mathcal{L}\psi T(\nabla \psi, \nabla u_{i})+ 4(T(\nabla \psi, \nabla u_{i}))^2\Big) \mathrm{dm}\\
     &= 4\|T(\nabla \psi, \nabla u_{i})\|^2_{L^2} - \int_\Omega u_{i}^{2} (\mathcal{L}\psi)^2\mathrm{dm} - 2\int_\Omega u_{i}^{2}T(\nabla \psi, \nabla \mathcal{L}\psi)\mathrm{dm}.
\end{align*}
Therefore,
\begin{align*}
    &\sum_{i=1}^k(\lambda_{k+1}-\lambda_{i})^2 \int_{\Omega}T(\nabla \psi, \nabla \psi)u_{i}^{2}\mathrm{dm} \\&\leq\sum_{i=1}^k(\lambda_{k+1}-\lambda_{i})\Big(4\|T(\nabla \psi, \nabla u_{i})\|^2_{L^2} - \int_\Omega u_{i}^{2} (\mathcal{L}\psi)^2\mathrm{dm} - 2\int_\Omega u_{i}^{2}T(\nabla \psi, \nabla \mathcal{L}\psi)\mathrm{dm}\Big).
\end{align*}
From \eqref{specialfunction-3} we know that $|\nabla \psi|^2=1$ and $\mathcal{L}\psi = B_0$, then using the previous inequality we obtain
\begin{align*}
    \varepsilon \sum_{i=1}^k(\lambda_{k+1}-\lambda_{i})^2 &\leq \sum_{i=1}^k(\lambda_{k+1}-\lambda_{i})(4\|T(\nabla u_{i})\|_{L^2}^2 - B_0^2)\\
     &\leq \sum_{i=1}^k(\lambda_{k+1}-\lambda_{i})(4\delta \lambda_{i} - B_0^2),
\end{align*}
which completes the proof of item ii).

iii) Taking $f=f_{\ell}$  into Lemma~\ref{lemma1} and summing over $\ell$, we get
\begin{align}\label{Equation-4.21}
   \sum_{i=1}^k(\lambda_{k+1}-\lambda_{i})^2& \int_\Omega \sum_{\ell=1}^{m+1}T(\nabla f_{\ell},\nabla f_{\ell})u_i^2\mathrm{dm} \nonumber\\
   &\leq\sum_{i=1}^k(\lambda_{k+1}-\lambda_{i}) \sum_{\ell=1}^{m+1} \int_\Omega (u_{i} \mathcal{L}f_{\ell}+2T(\nabla f_{\ell}, \nabla u_{i}))^2\mathrm{dm}.
    %& = \sum_{i=1}^k(\lambda_{k+1}-\lambda_{i}) \sum_{\ell=1}^{m+1}\int_\Omega (u_{i}^{2}(\mathcal{L}f_{\ell})^2 + 4u_{i}\mathcal{L}f_{\ell} T(\nabla f_{\ell}, \nabla u_{i}) + 4T(\nabla f_{\ell}, \nabla u_{i})^2) \mathrm{dm}
\end{align}
Using \eqref{specialfunction-2} and \eqref{property1} we have
\begin{equation}\label{Equation-(4.21)}
    \sum_{\ell=1}^{m+1}f_{\ell} \nabla f_{\ell} = 0, \quad \mbox{and} \quad \sum_{\ell=1}^{m+1}T(\nabla f_{\ell},\nabla f_{\ell}) = \gamma.
\end{equation}
Hence, from \eqref{Equation-(4.21)} and the Schwarz inequality, we obtain
\begin{align}\label{Eq_3.15}
    \gamma\sum_{i=1}^k(\lambda_{k+1}-\lambda_{i})^2 = \sum_{i=1}^k(\lambda_{k+1}-\lambda_{i})^2 \int_\Omega \sum_{\ell=1}^{m+1}T(\nabla f_{\ell}, \nabla f_{\ell})u_i^2\mathrm{dm},
\end{align}
\begin{align}\label{Equation-4.22}
    \sum_{\ell=1}^{m+1}\mathcal{L}f_{\ell} T(\nabla f_{\ell}, \nabla u_{i})= \sum_{\ell=1}^{m+1}-\gamma f_{\ell} T(\nabla f_{\ell}, \nabla u_{i})= -\gamma T\left(\sum_{\ell=1}^{m+1} f_{\ell} \nabla f_{\ell}, \nabla u_{i}\right)=0,
\end{align}
\begin{align}\label{Equation-4.23}
    \sum_{\ell=1}^{m+1} (\mathcal{L}f_{\ell})^2= \sum_{\ell=1}^{m+1}\big(-\gamma f_{\ell}\big)^2 = \gamma^2\sum_{\ell=1}^{m+1}f_{\ell}^2=\gamma^2,
\end{align}
and
\begin{align}\label{Equation-4.24}
    \nonumber
    \sum_{\ell=1}^{m+1} T(\nabla f_{\ell}, \nabla u_{i})^2&=\sum_{\ell=1}^{m+1} \langle \nabla f_{\ell}, T(\nabla u_{i})\rangle^2\leq\sum_{\ell=1}^{m+1} |\nabla f_{\ell}|^2|T(\nabla u_{i})|^2 \\
    \nonumber
    &\leq \sum_{\ell=1}^{m+1} \frac{1}{\varepsilon}\langle \nabla f_{\ell}, T(\nabla f_{\ell})\rangle\delta\langle \nabla u_{i}, T(\nabla u_{i})\rangle \\
    \nonumber
    &=\frac{\delta}{\varepsilon}\sum_{\ell=1}^{m+1}T(\nabla f_{\ell},\nabla f_{\ell})T(\nabla u_{i},\nabla u_{i}) \\
    &= \frac{\delta}{\varepsilon}\gamma T(\nabla u_{i},\nabla u_{i}).
\end{align}
Substituting \eqref{Eq_3.15}-\eqref{Equation-4.24} into \eqref{Equation-4.21}, we get
\begin{align*}
    \gamma\sum_{i=1}^k(\lambda_{k+1}-\lambda_{i})^2 &\leq \sum_{i=1}^k(\lambda_{k+1}-\lambda_{i})\Big(\gamma^2+4\frac{\gamma\delta}{\varepsilon}\int_{\Omega}T(\nabla u_{i},\nabla u_{i})\mathrm{dm}\Big)\\
    &\leq \gamma\sum_{i=1}^k(\lambda_{k+1}-\lambda_{i})\Big(\gamma+\frac{4\delta}{\varepsilon}\lambda_{i}\Big).
\end{align*}
Thus, the above inequality is enough to complete the proof of item iii) and finish the proof of Theorem~\eqref{theorem3}.
\end{proof}

\section*{Acknowledgements} 
The first author has been partially supported by Fundação de Amparo à Pesquisa
do Estado do Amazonas (FAPEAM).

\end{document}